\documentstyle{article}
\input{amssym.def}
\input{amssym.tex}
\begin{document}
%*******************************************************
\title{On volumes of hyperbolic 3-manifolds}
%*******************************************************

\author{Igor  ~Nikolaev\\
Department of Mathematics\\
2500 University Drive N.W.\\
Calgary T2N 1N4 Canada\\
{\sf nikolaev@math.ucalgary.ca}}

%**************************************************

%\date{January  12, 2003}
 \maketitle

%**************************************************

\newtheorem{thm}{Theorem}
\newtheorem{lem}{Lemma}
\newtheorem{dfn}{Definition}
\newtheorem{rmk}{Remark}
\newtheorem{cor}{Corollary}
\newtheorem{prp}{Proposition}
\newtheorem{exm}{Example}

%*************************************************

\newcommand{\N}{{\Bbb N}}
\newcommand{\F}{{\cal F}}
\newcommand{\R}{{\Bbb R}}
\newcommand{\Z}{{\Bbb Z}}
\newcommand{\C}{{\Bbb C}}

%******************************************************************
\begin{abstract}
The main thrust of present note is a volume formula
for hyperbolic surface bundle with the fundamental group $G$.
 The novelty consists in a purely algebraic
approach to the above problem. Initially, we concentrate
on the Baum-Connes morphism $\mu^G:K_{\bullet}(BG)\to K_{\bullet}(C^*_{red}G)$
for our class of manifolds, and then classify $\mu^G$ in terms 
of the ideals  in the ring of integers of a quadratic
number field $K$. Next, we extract the topological data
(e.g. volume of the orbifold ${\Bbb H}^3/G$) from the arithmetic
of field $K$.

\vspace{7mm}


\vspace{5mm}
{\it AMS (MOS) Subj. Class.:  46L85, 57M50}
\end{abstract}

%**************************************************************************
\section*{Introduction}
%***************************************************************************
The topological classification of $3$-dimensional manifolds is a challenge
of enormous interest and difficulty. Thanks to knot theory and
topological quantum field theory we understand better the nature 
of such manifolds. The lack of good invariants is usually blamed
for a slow progress in the field. Main expectations are connected
with Thurston's idea to substitute topology by geometry,
since most of the $3$-manifolds admit a hyperbolic metric.

It is amazing that geometry of 3-dimensional manifolds appears to be
essentially ``asymptotic'' or ``noncommutative''. Let us 
explain what  we mean by this. Heegaard splitting of a $3$-manifold $M$ 
uncovers one important detail about geometry of $M$. There
exist infinitely many (a countable set of) distinct Heegaard diagrams
which represent $M$. Every such diagram is a closed line on  
surface of genus $g\ge1$, so that the set of equivalent Heegaard
diagrams can be seen as a ``spectrum'' of $M$. It is clear
that any reasonable topological invariant of $M$ must depend exclusively
on the infinite, i.e. asymptotic, part of the spectrum. Note
that for $2$-dimensional manifolds (surfaces) this construction
always yields finite spectra.

Let $G=\pi_1M$, $BG$ the classifying space of $G$ and $C_{red}^*G$
the reduced group $C^*$-algebra of $G$. It was conjectured by Novikov, 
Baum and Connes that the assembly map
%**************************************************************
\begin{equation}\label{eq1} 
\mu^G:K_{\bullet}(BG)\longrightarrow K_{\bullet}(C_{red}^*G)
\end{equation}
%***************************************************************
is a morphism. Thus any $K$-theory invariant of 
$C^*_{red}G$ is an index, i.e. topological invariant of $M$.

Our knowledge of morphism (\ref{eq1}) is remarkably scarce.
Amenable groups, Gromov's hyperbolic groups, and some other
special groups are known to satisfy the Baum-Connes conjecture.
The general problem seems to be extremely hard. In this context,
the aim of present note is to calculate the Baum-Connes morphism
in case $G=\pi_1M$, where $M$ is a 3-dimensional manifold
fibering over the circle with surface $S$ as a fiber
(surface bundle).

This class of manifolds has been studied by Stallings, J\o rgensen
and Thurston. Stallings described the structure of fundamental
group of surface bundles. \linebreak 
J\o rgensen showed that for one punctured torus
bundle with aperiodic (pseudo-Anosov) monodromy, there exists a 
hyperbolic metric on $M$. Thurston proved that this is true for
every surface bundle with the pseudo-Anosov monodromy. On the other hand,
surface bundles seem to finitely cover {\it any} hyperbolic 3-manifold.

Let $G$ be the Kleinian group of surface bundle $M$. The assembly 
map $\mu^G$ admits in this case an explicit construction because
of the highly elaborated geometric theory of ``quasi-Fuchsian
deformations'' of the group $\pi_1S$ (Thurston and others). 
Roughly speaking, representations of $G$ and topology of $M$
are ``controlled'' by ``geometry'' of surface $S$. More precisely,
certain collections of disjoint curves on $S$, called {\it geodesic
lamination}, are responsible for the limit representation of
$\pi_1S$ in $SL(2,{\Bbb C})$. In turn, this limit representation 
generates the Kleinian group $G$ according to Stallings' result
for $\pi_1M$ (see Appendix). On the other hand, geodesic laminations
are intimately connected with the combinatorics of $AF$ $C^*$-algebras
 and through it with $C^*_{red}G$.

It can be shown by a routine calculation that $K$-theory
of $C^*_{red}G$ is Bott 2-periodic. We shall focus on the topological
data stored in the ordered $K_0$-group ${\cal E}=(K_0,K_0^+,[u])$,
known as {\it dimension group} of $C^*_{red}G$. For the surface
bundles, ${\cal E}$ has a remarkable property of being {\it
stationary}, i.e. self-similar
%***************************************************************
\footnote{In the context of holomorphic dynamics, this property
of surface bundles has been singled out by McMullen
as a ``renormalization'' property, see (\cite{McM}).}
%****************************************************************  
for an automorphism $A$ of lattice ${\Bbb Z}^n$.  The isomorphism class of 
stationary dimension group is known to intrinsically
depend on the ring $End~{\cal E}$ of the order-endomorphisms of $\cal E$. By a 
crucial calculation (Lemma \ref{lm2}) we show that
$End~{\cal E}\simeq O_K$, where $O_K$ is the ring of integers in a real
quadratic number field $K={\Bbb Q}(\sqrt{d})$. Thus, $O_K$ is a
Morita invariant of $K_0(C^*_{red}G)$.

We wish to explain the nature of ring $O_K$ in more geometric terms.      
Geodesic lamination $\Lambda$ comes with definite ``slope''
$\theta$ on surface $S$. Real number $\theta$ equals to average
inclination of the geodesic leaves in $\Lambda$, and can be
thought of as extension of the well-known Poincar\'e rotation
numbers to the case $S$ having genus $g\ge 2$ (see Appendix).
Since $\Lambda$ is fixed by a pseudo-Anosov map on $S$, its $\theta$
unfolds into a periodic continued fraction, and therefore
$\theta\in K={\Bbb Q}(\sqrt{d})$ for a positive square-free integer
$d$. In fact, $\theta\in O_K$ and all {\it cyclic covers} $\Lambda'$
of $\Lambda$, i.e. $\Lambda'$ such that $\psi(\Lambda')=\Lambda'$,
$\psi(\Lambda)\ne \Lambda$ and $\psi^m(\Lambda)=\Lambda$, will have
slopes $\theta'\in O_K$. Thus, the ring $O_K$ is generated by  the 
slopes of cyclic covers of lamination $\Lambda$.

\medskip\noindent
\underline{Volumes of hyperbolic 3-manifolds.}
The volume of surface bundle $M$ is the homotopy invariant
of $M$. Integration of hyperbolic
volumes consists in  an ``ideal triangulation'' of $M$,
so that the sides of ideal tetrahedra satisfy certain 
``compatibility equations''. It is remarkable that compatibility
equations are intrinsically diophantine and have clear links
to number theory (Neumann-Zagier \cite{NeZ}). Other methods
include numerical integration (SNAPPEA, Weeks \cite{Wee}),
and conjectures about coloured Jones polynomials (Kashaev \cite{Kas}) and Mahler 
measures (Boyd \cite{Boy}). An estimate of hyperbolic volumes through the 
Weil-Petersson metric can be given (Brock \cite{Bro}). 
For a class of ``arithmetic 3-manifolds'' the classical
Humbert's volume formula is known (Humbert \cite{Hum}).
Note that our approach indicates that surface bundles 
are ``arithmetic'', although it was  proved that 
there are only a finite number of surface bundles among
``classical'' arithmetic manifolds (Bowditch-Maclachlan-Reid
\cite{BMR}).    In general,  the problem of intergration seems to be 
 extremely hard.

\medskip\noindent
\underline{Scope of present note.}
In this note we introduce an integration technique for hyperbolic
surface bundles, inspired by the ideas and methods of noncommutative
geometry. The integration runs as follows. Let $N=N_{(p_1,q_1),\dots,(p_n,q_n)}$
be hyperbolic surface bundle obtained from a manifold $M$
with $n$ cusps by $(p_i,q_i)$-surgery on the $i$-th cusp.  
For brevity, consider the coprime pairs $p_i>0, q_i=1$
and take an infinite periodic continued fraction:
%************************************************************************
\begin{equation}\label{e12}
\theta=p_1+{1\over\displaystyle p_2 +\dots
{\strut 1\over\displaystyle p_n\displaystyle +\dots}}=Per~[p_1,\dots,p_n].
\end{equation}
%*****************************************************************  
Eventually, it can be proved that $\theta$ equals to the slope of geodesic
lamination $\Lambda$ on the fibre $S$ of surface bundle $N$. Take   
the quadratic extension $K=K(\theta)$ of field ${\Bbb Q}$. 
Finally, the last step of the integration process is
contained in the following theorem.

%************************************************************************
\smallskip\noindent
{\bf Theorem}
{\it 
Let $M$ be a hyperbolic surface bundle with $n$ cusps. Then the volume of
manifold $N=N_{(p_1,1),\dots,(p_n,1)}$ obtained from $M$ by the
$(p_i,1)$-filling of the $i$-th cusp is given by the formula:
%***************************************************************************
\begin{equation}\label{eq3*}
Vol~(N)=C(M){\log\varepsilon\over\sqrt{d}},
\end{equation}
%*****************************************************************
where $\varepsilon>1$ is the fundamental unit, $d$ the discriminant of
the field $K=K(p_1,\dots,p_n)$ and $C(M)>0$ a real constant depending
only on manifold $M$.   
}

%**********************************************************************
\bigskip\noindent
Let $\zeta_K(s)$ be the Dedekind zeta-function of field $K$. Then $\zeta_K(s)$
has a pole in the point $s=1$. It is not hard to see that formula
(\ref{eq3*}) can be written as $Vol~(N)=C(M) Res~\zeta_K(1)$ in this case. 
The latter expression is an analogue of Humbert's formula 
$Vol~(N)={|d|^{3\over 2}\over 4\pi^2}~~\zeta_K(2)$ for the volume of 
``arithmetic'' 3-manifold $N$, attached to the quadratic field $K$
with discriminant $d<0$.

%*********************************************************
%\tableofcontents
%**************************************************************************

%**************************************************************************
\section{Basic lemmas}
%***************************************************************************
This section contains calculation of the group $K_0(C^*_{red}G)$ (Lemma \ref{lm1})
and its invariants (Lemma \ref{lm2}). For the sake of brevity, we assume
that the reader is familiar with the 3-manifolds, assembly maps, $K$-theory of AF $C^*$-algebras
and algebraic number theory. For the reference we would recommend the
unpublished lecture notes of Thurston (\cite{T})
survey article of Higson (\cite{Hig}), books of R\o rdam {\it et al} (\cite{RLL}) and
Hecke (\cite{H}). Some very brief review can be found 
in the Appendix to this article. 
Let us introduce the following notation:

\bigskip
\begin{tabular}{cl}
${\Bbb H}^3$        & hyperbolic 3-dimensional space;\cr
                    &\cr
$G$                 & discrete (Kleinian) group on ${\Bbb H}^3$;\cr
                    &\cr
${\cal A}$          & $AF$ $C^*$-algebra;\cr
                    &\cr 
$\hat{\cal A}$      & representation space of algebra $\cal A$;\cr
                    &\cr
${\cal E}$          & dimension (Elliott) group;\cr
                    &\cr
$\theta_{\cal E}$   & rotation number associated to ${\cal E}$;\cr
                    &\cr
${\cal E}_{A}$      & stationary dimension group given by  matrix $A$;\cr
                    &\cr
$B({\cal A})$       & Bratteli diagram of ${\cal A}$;\cr
                    &\cr
$K_0$               & Elliott functor ${\cal A}\to {\cal E}$;\cr
                    &\cr
$\Lambda$           & geodesic lamination on surface $S$;\cr
                    &\cr
$C^*_{red}G$              & reduced group $C^*$-algebra.
\end{tabular}

%***********************************************************************
\subsection{Calculation of the group $K_0(C^*_{red}G)$}
%*************************************************************************
\begin{lem}\label{lm1}
If $G$ is the fundamental group of  
3-dimensional manifold fibering over the circle with a pseudo-Anosov 
monodromy $\varphi$, then
%*****************************************************************
\begin{equation}\label{eq2}
K_0(C^*_{red}G)\cong{\cal E}_A.
\end{equation}
%********************************************************************
\end{lem}
%*************************************************************************
{\it Proof.} We split the proof into several parts (Propositions
\ref{prp1}-\ref{prp3}), which might be of independent interest 
to the reader. We apologize for a fairly sketchy argument,
the reader can find details in the original articles and
preprints.
%*******************************************************************
\begin{prp}\label{prp1}
{\bf (Bers-Thurston)}
Given pseudo-Anosov diffeomorphism $\varphi: S\to S$,  
the Kleinian group $G$ of Lemma \ref{lm1} is determined by a
geodesic lamination $\Lambda\subset S$ such that 
%***************************************************
\begin{equation}\label{eq3}
\varphi(\Lambda)=\Lambda.
\end{equation}
%*****************************************************
\end{prp}
%*****************************************************************
{\it Proof.} The original proof of this fact is due to
Thurston (\cite{Thu}). Let $M$ be the mapping torus of a surface
homeomorphism $\varphi:S\to S$. Thurston proved (\cite{Thu}, Theorem 5.5)
that $M$ is either a Seifert fibration ($\varphi$ is periodic) or splits
into simpler manifolds ($\varphi$ is reducible)  or else a hyperbolic
3-manifold ($\varphi$ is pseudo-Anosov). In the latter case, there
exists a Kleinian group $G$, which is a ``geometric limit'' of
quasi-fuchsian deformations of the group representations
$\pi_1S\to SL(2,{\Bbb C})$. This limit has a remarkable description
in terms of a $\varphi$-invariant geodesic lamination $\Lambda$
on surface $S$. (In other words, geometry on $S$ ``controls''
topology of $M$.) By the Mostow-Prasad rigidity of group
representations, $\Lambda$'s are in a 1-1 correspondence with
groups $G$ (\cite{Thu}). Since $G\cong \pi_1(S)\rtimes_{\varphi}
{\Bbb Z}$, Proposition \ref{prp1} follows.  
$\square$

%*******************************************************************
\begin{prp}\label{prp2}
(\cite{Nik})
For any $\cal A$ there exists $\Lambda$ such that
$\hat {\cal A}\cong\Lambda$. Conversely, given $\Lambda$
there exists a dual $\hat\Lambda\cong {\cal A}$.
\end{prp}
%*****************************************************************
{\it Proof.} The idea of proof is
based on bijection between Bratteli diagrams of $AF$ $C^*$-algebras and 
geodesic laminations on $S$, which follows, in turn, from the Koebe-Morse
theory of coding of the geodesic lines. 
Let $\cal A$ be a simple $AF$ $C^*$-algebra and $B({\cal A})$ its Bratteli
diagram of rank $2g$. Consider the set $Spec~B({\cal A})$ consisting
of all infinite paths of the diagram $B({\cal A})$. Marking the vertices
(of same level) of $B({\cal A})$ by the ``symbols'' 
${\goth A}=\{{\goth a}_1,\dots,{\goth a}_{2g}\}$, one
gets an infinite word
$\{\omega_1,\omega_2,\dots~|~ \omega_i\in {\goth A}\}$  
for every element of $Spec~B({\cal A})$. Since each word
 is a Morse symbolic geodesic,
one gets a correspondence between $Spec~B({\cal A})$ and a geodesic
lamination $\Lambda$ on surface $S$. Note that
$\Lambda$ is minimal since any two geodesics lie in the closure
of each other by the simplicity of ${\cal A}$. 

On the other hand, every infinite path of $B({\cal A})$ defines
a pure state on the $C^*$-algebra $\cal A$. In terms of the representation
theory,  this means that $\hat A\cong Spec~B({\cal A})$. Proposition
\ref{prp2} follows. 
$\square$

%*******************************************************************
\begin{prp}\label{prp3}
Let $\Lambda$ be geodesic lamination satisfying equation 
(\ref{eq3}) and let $\hat\Lambda\cong {\cal A}$
be its dual (Proposition \ref{prp2}). Then
%***************************************************
\begin{equation}\label{eq5}
K_0({\cal A})={\cal E}_A,
\end{equation}
%*****************************************************
for a stationary dimension group ${\cal E}_A$.
\end{prp}
%*****************************************************************
{\it Proof.} We have to show that any geodesic lamination on $S$,
which is invariant under a pseudo-Anosov diffeomorphism $\varphi$,
must be given by Bratteli diagram of stationary type.
Every Bratteli diagram
$B({\cal A})$ admits a natural order $\ge$ on its set of infinite
paths $Spec~B({\cal A})$, see (\cite{HPS}). 
Based on this order, one can define
a Putnam-Vershik homeomorphism on the set $Spec~B({\cal A})$, which is a Cantor set
in a natural topology. If the diagram is simple, the homeomorphism
is minimal.

Consider a restriction $\phi=\varphi~|_{\Lambda}$ of $\varphi$ to geodesic
lamination $\Lambda\cong Spec~B({\cal A})$. It is not hard to see
that $\phi$ is conjugate to a Putnam-Vershik homeomorphism
on the Cantor set $\Lambda$. Moreover, the condition $\varphi (\Lambda)
=\Lambda$ would imply that $\phi$ must be order-preserving with
respect to the ordered diagram $(B({\cal A}),\ge)$. The latter
requirement can be met if and only if the diagram admits 
strong structural ``symmetry'' which is actually a periodicity
with some constant incidence matrix $A$. 
Proposition \ref{prp3} follows.   
$\square$

\bigskip\noindent
To finish the proof of Lemma \ref{lm1}, let us notice that
$C^*G\simeq \hat G$, where $\hat G$ is a ``noncommutative dual''
of the Kleinian group $G$. Using Propositions \ref{prp1}
and \ref{prp2} we get the following chain of isomorphisms:
%******************************************************************
\begin{equation}\label{eq6}
C^*_{red}G \cong  \hat G \cong \hat\Lambda\cong {\cal A},
\end{equation}
%*****************************************************************
and $K_0(C^*_{red}G) = K_0({\cal A})={\cal E}_A$ by Proposition
\ref{prp3}. Lemma \ref{lm1} follows.
$\square$

%***********************************************************************
\subsection{Classification of groups $K_0(C^*_{red}G)$}
%*************************************************************************
Let ${\cal E}=(E,E^+,[u])$ be ordered abelian group of $rank~E<\infty$.
Consider the ring $End~E$ of endomorphisms of $E$ under ``pointwise''
addition and multiplication of mappings $E\to E$. The subring of
$End~E$, made of endomorphisms which ``preserve'' positive cone $E^+$,
we shall denote by $End~{\cal E}$. 

Let ${\cal E}_A$ be stationary ordered abelian group. In this case,
the rotation number $\theta_{{\cal E}_A}$ is known to be real
quadratic (see Appendix). Denote by $K={\Bbb Q}(\theta_{{\cal E}_A})$
the field of algebraic numbers, which is an extension of $\Bbb Q$
by quadratic irrationality $\theta_{{\cal E}_A}$. Clearly,
$K={\Bbb Q}(\sqrt{d})$ for a square-free positive integer $d$.

Finally, consider the embedding 
%******************************************************************
\begin{equation}\label{eq7}
End~{\cal E}_A\rightarrow K,
\end{equation}
%*****************************************************************
given by the formula ${\cal E}_i\mapsto \theta_{{\cal E}_i}$,
where ${\cal E}_i=Im~f_i$ and $f_i\in End~{\cal E}$. It is not
hard to see that full image of $End~{\cal E}_A$ is an integral
ideal in $K$, which we denote by $I_{{\cal E}_A}$. 
%***************************************************************
\begin{dfn}\label{dfn1}
The equivalence class
%**************************************************************
\footnote{Two integral ideals $I_1,I_2$ are equivalent if and only
if there exists a principal ideal $(\omega),\omega\in O_K$ such that
$I_2=(\omega)I_1$, see Hecke (\cite{H}).}
%*************************************************************
of ideals $[I_{{\cal E}_A}]$
in the ring of integers $O_K$ we call associated to
the stationary dimension group ${\cal E}_A$. 
\end{dfn}
%*************************************************************
%***********************************************************************
\begin{lem}\label{lm2}
Let $[{\cal E}_A]$ be the Morita equivalence class of stationary
dimension group ${\cal E}_A$. The correspondence $[{\cal E}_A]
\leftrightarrow [I_{{\cal E}_A}]$ is a bijection which classifies
stationary dimension groups.
\end{lem}
%***********************************************************************
{\it Proof.} 
The idea of the proof belongs to Handelman, cf. \S 5 of (\cite{Han}).
We split the proof into two propositions.
Since ${\cal E}_A$ is stationary, its ring 
of order-preserving endomorphisms $End~{\cal E}$ is non-trivial, i.e. distinct 
from ${\Bbb Z}$. Morover, in this case $End~{\cal E}_A$ completely
defines, and is defined by, ${\cal E}_A$ (Proposition \ref{prp4}). 
We pull back the embedding (\ref{eq7}) and prove that Morita equivalent
dimension groups  generate equivalent ideals in $O_K$ (Proposition \ref{prp5}).
%*******************************************************************
\begin{prp}\label{prp4}
Stationary dimension group ${\cal E}_A$ is defined by its endomorphism
ring $End~{\cal E}_A$.
\end{prp}
%*******************************************************************
{\it Proof.} Let ${\cal E}_A$ be stationary dimension group of rank
$n$. Consider the set
%******************************************************************
\begin{equation}\label{eq8}
X=\bigcap_{f\in End~{\cal E}_A} f({\cal E}_A).
\end{equation}
%******************************************************************
$X$ is a non-empty set, since $End~{\cal E}_A\ne {\Bbb Z}$. It is also
an additive abelian group of rank $n$. Clearly $f\in End~{\cal E}_A$
coincide with $f\in Aut~X$. It is not hard to see that $Aut~X$ is
infinite cyclic group generated by a single automorphism $f:{\Bbb Z}^n\to
{\Bbb Z}^n$. Then our stationary group ${\cal E}_A$ is representable as limit 
of the simplicial dimension groups under the automorphism $f$. 
$\square$

%*******************************************************************
\begin{prp}\label{prp5}
Let ${\cal E}_A$ and ${\cal E}_A'$ be two (stationary) dimension groups,
which are Morita equivalent. Then their associated ideals $I_{{\cal E}_A}$
and $I_{{\cal E}_A'}$ are equivalent, i.e. differ by a principal ideal
 multiple $(\omega)$.
\end{prp}
%*******************************************************************
{\it Proof.} 
Let ${\cal E}_A=({\Bbb Z}^n,{\Bbb Z}^n_+,[u])$ and
${\cal E}_A'=({\Bbb Z}^n,({\Bbb Z}^n_+)',[u'])$ be two
stationary dimension groups of rank $n$. Then 
${\cal E}_A\sim {\cal E}_A'$ are Morita equivalent if and
only if there exists $f\in End~{\cal E}_A$ such that 
$f({\Bbb Z}^n_+)=({\Bbb Z}^n_+)'$. Since every $f$ extends
to a linear map $\widetilde f:{\Bbb R}^n\to
{\Bbb R}^n$, one can consider the $n$-torus $T^n={\Bbb R}^n/{\Bbb Z}^n$
and the mapping:
%**************************************************************
\begin{equation}\label{eq9}
\widetilde f: T^n\longrightarrow T^n.
\end{equation}
%******************************************************************
Let now $I_{{\cal E}_A}$ be the ideal associated to ${\cal E}_A$.
Fix an algebraic integer $\omega\in O_K$ such that 
%**************************************************************
\begin{equation}\label{eq10}
N(\omega)=deg~\widetilde f,
\end{equation}
%******************************************************************
where $N(\omega)$ is the norm of $\omega$ and $deg~\widetilde f$ 
is the degree of continuos map (\ref{eq9}). Since the norm is 
a multiplicative function on the ideals and the norm of principal ideal $(\omega)$ 
coincides with the $N(\omega)$, one gets that  $I_{{\cal E}_A'}=(\omega)I_{{\cal E}_A}$.
$\square$

\bigskip\noindent
To finish the proof of Lemma \ref{lm2}, it remains to apply 
Propositions \ref{prp4} and \ref{prp5}, and recall the definition of
the equivalence class of the ideal.
$\square$

%**************************************************************************
\section{Main result}
%***************************************************************************
In this section we state main result which relates hyperbolic volumes of
surface bundles with the arithmetic of field $K$ (volume inequalities). 
Our proof of these inequalities  is based on an observation that ``commensurability
class'' of hyperbolic manifold has a natural ``ideal structure'', which is
isomorphic to the ideal structure  (``arithmetic'') of field $K$. 

Before we go to the statement of results, some general remarks on the analogy
between number theory and 3-dimensional topology might be hepful.
Hyperbolic surface bundles might be good example of ``noncommutative 
varieties''. Indeed, think of ${\cal E}_A$ as a ``coordinate ring''
of the noncommutative variety $V({\cal E}_A)$. Then ``divisor'' on
$V({\cal E}_A)$ coincides with an ideal in ring $O_K$ (Lemma \ref{lm2}).
The group of equivalence classes of divisors $Pic~V({\cal E}_A)=Cl~K$,
where $Cl~K$ is the class group of field $K$. The role of the Chern class
of the line bundle $L_{V({\cal E}_A)}$ is played by the ``volume formula''
$Ch~ L_{V({\cal E}_A)}\simeq Res~\zeta_K(1)$, where $\zeta_K$ is the Dedekind
zeta-function of field $K$.

%***************************************************************
\begin{dfn}\label{dfn2}
Manifolds $M_1$ and $M_2$ are said to be commensurable,
if there exist finite coverings $\widetilde M_1$ and $\widetilde M_2$, 
such that $\widetilde M_1\simeq\widetilde M_2$ are
homeomorphic. The equivalence class of manifolds commensurable with
$M$ is denoted by ${\goth M}$. Manifold $M\in {\goth M}$ is called prime, 
if it cannot be covered by any other manifold from ${\goth M}$ except 
$M$ itself.
\end{dfn}
%*************************************************************
Let ${\goth M}$ be commensurability class of a hyperbolic surface bundle.
Lemma \ref{lm2} indicates on an   isomorphism between ``arithmetic'' of 
set ${\goth M}$ and the ideal structure of ring $O_K$. This observation
is by no means new in the context of ``arithmetic'' 3-manifolds,
cf. Helling (\cite{Hel}). Borel (\cite{Bor}) showed that there exists 
infinitely many prime manifolds in every commensurability class
of such manifolds. Thurston observed that one can study invariants
 of  commensurability  classes of 3-dimensional manifolds using 
certain algebraic number fields (trace-fields), cf \S 6.7 of (\cite{T}). 
The following theorem is an attempt to extend Borel-Helling-Thurston
theory to hyperbolic surface bundles.
%***************************************************************
\begin{thm}\label{thm1}
Let $M$ be hyperbolic surface bundle, which has minimal volume
in its commensurability class ${\goth M}$. Then there exit real
constants $k$ and $K$ depending exclusively on ${\goth M}$, such
that
%**************************************************************
\begin{equation}\label{eq12}
k{\log\varepsilon\over\sqrt{d}}\le Vol~M\le K{\log\varepsilon\over\sqrt{d}},
\end{equation}
%******************************************************************
where $\varepsilon$ is the fundamental unit of the number field
$K={\Bbb Q}(\sqrt{d})$ associated to $M$. Moreover,
if $\{N_i\}$ is the family of surface bundles obtained from
$M$ by the Dehn surgery (cusp-filling), then in formula (\ref{eq12})
$k=K=C(M)>0$ for a positive constant $C(M)$ defined for the whole family
$\{N_i\}$ and depending only on manifold $M$.  
\end{thm}
%**************************************************************
%***************************************************************
\begin{cor}\label{cor1}
Function $M\to Vol~M$ has degree $h_K$ at $M$, where $h_K$ is the class
number of field $K$.
\end{cor}
%**************************************************************

%**************************************************************************
\section{Proof}
%***************************************************************************
%***********************************************************************
\subsection{Proof of Theorem 1}
%***********************************************************************
{\sf Part I.} We wish to prove inequality (\ref{eq12}).
Let us outline the main idea. If $M$ is a surface bundle with the
pseudo-Anosov monodromy, then it can be described as an ideal in
the equivalence class of ideals $[I_{{\cal E}_A}]$. Other ideals
in this class correspond to manifolds commensurable with $M$. Thus,
one can extend the notion of ``divisibility'' from the
ideals to commensurable manifolds. In this sense, prime ideals
correspond to prime manifolds in the commensurability class
${\goth M}$. By the Borel-Helling
lemma the number of prime manifolds is infinite (Lemma \ref{lm3}). 
Note also that since $O_K$ is the Dedekind domain, any $M\in {\goth M}$ can 
be ``uniquely decomposed into primes''. On the other hand, the ``volume
growth'' of manifolds in ${\goth M}$ is discrete, in sense that it
is bounded by the uniform constants $k$ and $K$ (Lemma \ref{lm4}).
Thus, knowing the Dirichlet density of ideals in $[I_{{\cal E}_A}]$
(Lemma \ref{lm5}) allows us to estimate the volume of minimal prime
manifold in ${\goth M}$ in terms of the arithmetic of field $K$.

\bigskip\noindent
\underline{Division algorithm for commensurable manifolds}. 
Let ${\goth M}$ be commensurability class of hyperbolic surface
bundles. If $M_1,M_2\in {\goth M}$ are such that $M_1$ covers
$M_2$, one can define their ``ratio'' $M_1/M_2\in {\goth M}$
in the following way. Let ${\cal E}_{A_1}$ and ${\cal E}_{A_2}$
be the stationary dimension groups corresponding to $M_1$ and
$M_2$ (Lemma \ref{lm1}). Since $M_1$ covers $M_2$, there exists
a order-preserving homomorphism $f:{\cal E}_{A_1}\to {\cal E}_{A_2}$.
Then the normal (order) subgroup ${\cal E}_{A_1'}=Ker~f\subseteq {\cal E}_{A_1}$ is 
a stationary dimension group. Thus if $M_1'$ is the
hyperbolic surface bundle corresponding to ${\cal E}_{A_1'}$, we can
define the ratio: 
%**************************************************************
\begin{equation}\label{eq13}
{M_1\over M_2}=M_1'.
\end{equation}
%******************************************************************
Clearly, $M_1'\in {\goth M}$ and either $M_1'$ covers $M_2$
or ``relatively prime'' with $M_2$. In the first case, we repeat
the division algorithm until the newly obtained ratio becomes
``relatively prime'' to $M_2$. By the construction, the
algorithm stops on the $p$-th step, if and only if $M_1$ is  
$p$-fold cover to $M_2$.

%***********************************************************************
\begin{lem}\label{lm3}
{\bf (Borel-Helling)}
The number of prime manifolds in the commensurability class ${\goth M}$
is infinite.
\end{lem}
%***********************************************************************
{\it Proof.} For hyperbolic manifolds, which are factor spaces of ${\Bbb H}^3$
by the discrete subgroup of $SL(2,{\Bbb C})$ generated by the algebraic
integers of the imaginary quadratic field $K={\Bbb Q}(\sqrt{-d})$, this
lemma was proved by Borel (\cite{Bor}) and Helling (\cite{Hel}). We wish
to prove it for hyperbolic surface bundles using the above mentioned
``division algorithm''. 

Let $M\in {\goth M}$ be a hyperbolic surface bundle of commensurability
class ${\goth M}$. By the division algorithm, $M$ unfolds into a product
of powers of prime manifolds: $M=(M_{\pi_1})^{p_1}(M_{\pi_2})^{p_2}\dots
(M_{\pi_k})^{p_k}$. Let $P_{\pi_i}$ be the ideals in $K$, which correspond
to prime manifolds $M_{\pi_i}$. It is not hard to see 
that $P_{\pi_i}$ are prime ideals in $K$. 

Let $\omega\in O_K$ be prime algebraic integer which does not belong
to any of the ideals $P_{\pi_i}, i=1,\dots,p$. Consider the principal
ideal $(\omega)$ and take the manifold $M_{\omega}\in {\goth M}$, which 
correspond to $(\omega)$. Clearly, $M_{\omega}$ is prime manifold,
which does not belong to finite list of prime manifolds 
$M_{\pi_1},M_{\pi_2},\dots, M_{\pi_k}$. 
$\square$

%***********************************************************************
\begin{lem}\label{lm4}
Let $M_0<M_1<M_2<\dots$ be sequence of commensurable hyperbolic 3-manifolds
of growing volume. Then there exists positive constants $k$ and $K$ such that for
every $i=1,\dots,\infty$ it holds:
%**************************************************************
\begin{equation}\label{eq14}
k\le Vol~M_i- Vol~M_{i-1}\le K.
\end{equation}
%****************************************************************** 
\end{lem}
%***********************************************************************
{\it Proof.} {\sf Step 1. Lower bound.}
To prove existence of $k>0$, it is sufficient to show that 
monotone sequence $Vol~M_0,Vol~M_1,\dots$ cannot be a Cauchy
sequence. (Indeed, for otherwise, $k=0$.) Suppose to the contrary, 
that monotone growing sequence $Vol~M_0,Vol~M_1,\dots$ is Cauchy.
Then volume function $Vol: {\goth M}\to {\Bbb R}$ is bounded by a
constant $C$. This gives us a contradiction, since any $M\in {\goth M}$
has an $m$-fold covering, such that $m~Vol~M>C$ for $m$ sufficiently 
large.

\medskip
{\sf Step 2. Upper bound.}
Let $K$ be a constant, such that $Vol~M_0< K$. Suppose to
the contrary, that there exist $i$ so that $Vol~M_{i+1}-
Vol~M_i>K$. Let $m$ be the maximal integer which satisfy
the equality $m Vol~M_0 \le Vol~M_i$. (In other words,
one takes the maximal cover of $M_0$, whose volume does
not exceed volume of $M_i$.) It is not hard to see that
in view of condition $Vol~M_0< K$ we have
%**************************************************************
\begin{equation}\label{eq15}
Vol~M_i< (m+1) Vol~M_0< Vol~M_{i+1}.
\end{equation}
%******************************************************************
Thus, $(m+1)$ covering of $M_0$ lies ``between'' $M_i$ and $M_{i+1}$,
what contradicts our assumptions. Lemma \ref{lm4} is proved.
$\square$

%***********************************************************************
\begin{lem}\label{lm5}
{\bf (Dirichlet density)}
Let $A$ be an equivalence class of ideals in the real
quadratic number field $K$ with the discriminant $d$ and fundamental
unit $\varepsilon$. Denote by $N(t, A)$ the number
of ideals ${\goth a}\in A$ such that $N{\goth a}\le t$, 
where $t$ is a positive integer. Then
%************************************************************
\begin{equation}\label{eq6.55}
\lim_{t\to\infty}{N(t, A)\over t}= {2\log\varepsilon\over\sqrt{d}}.
\end{equation}
%*************************************************************
Moreover, the above limit exists and is the same for all equivalence
classes of ideals in $K$.
\end{lem}
%****************************************************************   
{\it Proof.} For the proof of this lemma we refer the reader to
Hecke (\cite{H}).
$\square$

\bigskip
Let us complete the proof of Theorem \ref{thm1}. 
Denote by $N(t)=N(t,A)$ the number of ideals in the ideal
class $A$ whose norm doesn't exceed $t$ and consider a
chain of manifolds $\{M_0,M_1,\dots,M_{N(t)}\}\in {\goth M}$ 
of growing volume. The telescoping sequence
%**************************************************************
\begin{eqnarray}\label{eq17}
Vol~M_{N(t)} &-& Vol~M_{N(t)-1}+Vol~M_{N(t)-1}-Vol~M_{N(t)-2}+\dots
\cr\nonumber
+Vol~M_2 &-& Vol~M_1+Vol~M_1-Vol~M_0
 = M_{N(t)}-M_0,
\end{eqnarray}
%******************************************************************
can be evaluated using Lemma \ref{lm4} as:
%**************************************************************
\begin{equation}\label{eq18}
k(N(t)-1)< Vol~M_{N(t)}- Vol~M_0< K(N(t)-1).
\end{equation}
%******************************************************************
On the other hand, one has a prime decomposition of manifold
$M_{N(t)}$ in its commensurability class $\goth M$:
%**************************************************************
\begin{equation}\label{eq19}
M_{N(t)}=(M_{\pi_1})^{p_1}(M_{\pi_2})^{p_2}\dots(M_{\pi_k})^{p_k},
\end{equation}
%******************************************************************
where $M_{\pi_1}=M_0$ in view of Lemma \ref{lm3}.
Since we didn't make any assumption on the positive integer $t$
so far, we can now make one. Namely, take $t$ be such that $t=p_1$.   
One can now evaluate that
%**************************************************************
\begin{equation}\label{eq20}
Vol~M_{N(t)} =t Vol~M_0.
\end{equation}
%******************************************************************
Together with equation (\ref{eq18}), the latter will give us
that
%**************************************************************
\begin{equation}\label{eq21}
k(N(t)-1)< (t-1) Vol~M_0< K(N(t)-1).
\end{equation}
%******************************************************************
Since
%**************************************************************
\begin{equation}\label{eq22}
\lim_{t\to\infty}{N(t)-1\over t-1}=\lim_{t\to\infty}{N(t)\over t}=
{2\log\varepsilon\over\sqrt{d}},
\end{equation}
%******************************************************************
we get inequality (\ref{eq12}). Part I of Theorem \ref{thm1} is proved.

\bigskip\noindent
{\sf Part II.} By Thurston's result, we have $Vol~(N_i)\rightarrow Vol~(M)$
and $Vol~(N_i)<Vol~(M)$ as $i\to\infty$ by the cusp-filling process, 
see e.g. (\cite{NeZ},\cite{T}). On the other hand, constants $k_i$ and $K_i$
for $N_i$ in (\ref{eq12}) are exact and therefore $k_i=K_i=C(N_i)$ for $N_i$ whose
volume is infinitely close to $Vol~(M)$. Thus, $Vol~(N_i)=C(N_i)(\log\varepsilon_i/
\sqrt{d_i})$. 

In fact, $C(N_i)=Const=C(M)$ for the whole family $\{N_i\}$. To see
this, one should take a sequence of discriminants $\{d_i\}$,
such that the Dirichlet density of the fields $K_i={\Bbb Q}(\sqrt{d_i})$
is monotone growing and bounded. Since there is a subsequence in
$\{N_i\}$ corresponding to $\{K_i\}$, then by Thurston's formula
$Vol~(N_i)\rightarrow Vol~(M)$, we get $C(N_i)=Const$. Part II of 
 Theorem \ref{thm1} is proved.
$\square$

%***********************************************************************
\subsection{Proof of Corollary 1}
%***********************************************************************
For the imaginary quadratic number fields this fact was established
by Bianchi (\cite{Bia}). To establish similar result for
real quadratic number fields, let us notice that the Dirichlet density
does not depend on the choice of the equivalence class of ideals in
$O_K$ (Lemma \ref{lm5}). On the other hand, we have the following lemma.
%***********************************************************************
\begin{lem}\label{lm6}
Let $[{\cal E}_A]$ be the Morita equivalence class of stationary
dimension group ${\cal E}_A$. Let $K={\Bbb Q}(\sqrt{d})$ be the
associated number field to  ${\cal E}_A$. If the class number 
$h_K$ of the field $K$ is bigger than $1$, there exits $h_K$
distinct Morita equivalence classes of stationary dimension
groups with the same associated number field $K$. 
\end{lem}
%************************************************************
{\it Proof.}
Since $h_K$ is the number of equivalence classes of ideals in
the ring of integers $O_K$, one can apply Lemma \ref{lm2}. 
Then there exists exactly $h_K$ stationary dimension groups
${\cal E}_A^{(1)},{\cal E}_A^{(2)},\dots,{\cal E}_A^{(h_K)}$,
which are pairwise Morita non-isomorphic, but have the same 
associated number field $K$. Therefore, there exists the 
same number of Morita equivalence classes related to the 
given field $K$.  
$\square$

\medskip
To finish the proof of corollary, let $M\to Vol M$ be the
Gromov-Thurston function on the set of hyperbolic surface
bundles. From the construction of $M$ it follows that
stationary dimension groups ${\cal E}_A^{(1)},{\cal E}_A^{(2)},\dots,{\cal E}_A^{(h_K)}$
generate surface bundles $M^{(1)},M^{(2)},\dots,M^{(h_K)}$, which are topologically
distinct. However, the arithmetic of their commensurability classes is the
same for all of these manifolds. By the discussion above, the arithmetic
determines volume of the hyperbolic manifolds. Corollary \ref{cor1} follows.
$\square$

%***********************************************************************
\section{Appendix}
%**********************************************************************
%***********************************************************************
\subsection{Kleinian groups and hyperbolic 3-manifolds}
%***********************************************************************
We briefly review Thurston's theory of Kleinian groups for surface
bundles with the pseudo-Anosov monodromy. Excellent source
of original information is Thurston's paper (\cite{Thu}).

\medskip\noindent
\underline{Mapping tori.} Let $\varphi$ be a diffeomorphism
of surface $S$. One can produce a 3-manifold $M_{\varphi}$ depending 
on diffeomorphism $\varphi$ by identification of points
$(x,1)$ and $(\varphi(x),0)$ in the product $S\times I$ ($I$ - unit
interval). By Thurston's theory of surface diffeomorphisms (\cite{Thu}),
any $\varphi$ falls into one of three classes: (i) periodic, (ii) non-periodic
(pseudo-Anosov) and (iii) reducible (combination of two preceding
types). We shall be interested in $M_{\varphi}$, where $\varphi$ is
pseudo-Anosov. Such 3-manifolds are known to be fibre bundles over the circle,
with fibre $S$, and pseudo-Anosov monodromy $\varphi$.   
The interest to this case lies in the fact that $M=M_{\varphi}$ can
be presented as $M={\Bbb H}^3/G$, where $G=\pi_1M$ is a Kleinian
group. In other words, $M$ is a hyperbolic manifold of finite volume.
In fact, the topology of $M$ and representation theory of group $G$  
is perfectly ``controlled'' by geometry of surface $S$, hidden in the
notion of ``geodesic lamination'' on $S$.

\medskip\noindent
\underline{Geodesic laminations.}
Let $S={\Bbb H}^2/\Gamma$ be a hyperbolic surface 
given by the action of Fuchsian group $\Gamma$ 
on the hyperbolic plane ${\Bbb H}^2$. By a geodesic
line on $S$ one understands a line consisting
of geodesic (locally shortests) arcs in the given
hyperbolic metric on $S$. The geodesic line
is called {\it simple} if it has no transversal self-intersection
on $S$. {\it Geodesic lamination} is a disjoint union of
simple geodesic lines on surface $S$. 
If $\gamma$ is non-periodic geodesic,  then the topological
closure $Clos~\gamma$ on $S$ , contains  continuum of
non-periodic  geodesic lines,  having the same closure as $\gamma$. 
The set of all simple non-periodic geodesic lines that  are
everywhere dense in the set $Clos~\gamma$ we call {\it minimal
geodesic lamination} on $S$.

\medskip\noindent
\underline{Representation theory of Fuchsian groups.}
Let $\Gamma$ be Fuchsian (surface) group. Let us
consider the set $Rep~\Gamma$ of all faithful representations $\rho:\Gamma\to
SL_2({\Bbb C})$ into the group of isometries of hyperbolic
space ${\Bbb H}^3$. According to Thurston's theory
for $Rep~\Gamma$, up to conjugacy by isometries of ${\Bbb H}^3$
there exist three kinds of $\rho(\Gamma)$: (i) Fuchsian, i.e. 
when $\rho(\Gamma)$ preserves the boundary of unit disc $D\subset\partial {\Bbb H}^3$;
(ii) quasi-Fuchsian, i.e. when $\rho(\Gamma)$ ``deforms'' $D$ but
its boundary $\partial D$ is still a Jordan curve (topological circle)
in $\partial {\Bbb H}^3$; and (iii) discontinuous, i.e. when
$\rho(\Gamma)$ ``breaks'' the boundary $\partial D$ to a plane-filling
curve, homeomorphic to ${\Bbb R}$. The invariants of the type (i) and
(ii) representations are so-called (pair of) {\it conformal structures} 
on $S$, while the main invariant of type (iii) representation is an
{\it ending (geodesic) lamination} on $S$.

\medskip\noindent
\underline{Connection to hyperbolic geometry of mapping tori.}
The connection of $Rep~\Gamma$ to the mapping tori with pseudo-Anosov
monodromy is based on the Stallings' theorem about  structure
of the fundamental group of manifolds which fibre over the circle. 
Namely, the fundamental group of $M_{\varphi}$
($\varphi$ pseudo-Anosov) has a representation:  
%**************************************************************
\begin{equation}\label{eq23}
\pi_1M_{\varphi}=\langle \pi_1S, t ~|~tgt^{-1}=\varphi_*(g),~\forall g\in\pi_1S\rangle, 
\end{equation}
%******************************************************************
where $\varphi_*$ is the action of $\varphi$ on $\pi_1S$.Thus, to
construct a representation $\pi_1M_{\varphi}\to PSL_2({\Bbb C})$
(Kleinian group), is equivalent to construction of such a 
$\rho^*: \pi_1S\to PSL_2({\Bbb C})$, which is ``invariant''
(up to an isometry of ${\Bbb H}^3$) under the action of $\varphi_*$.  
Note that there is no chance $\rho^*$ to be representation of
type (i) or (ii). In fact, it is a basic observation of W.~P.~Thurston that 
$\rho^*$ must be of type (iii) (Thurston \cite{Thu}). Thus, ending 
geodesic laminations on $S$ ``controls'' topology of manifold $M_{\varphi}$. 
The following general result is true.    
%**************************************************************
\begin{thm}
{\bf (W.~P.~Thurston)}
Let $M$ be surface bundle with a pseudo-Anosov monodromy 
$\varphi: S\to S$. Then $M$ is a hyperbolic 3-manifold,
whose Kleinian group $G\cong \pi_1M$ is presented by a
$\varphi$-invariant (measured) geodesic lamination on surface $S$.  
\end{thm}
%*************************************************************
{\it Proof.} See (\cite{Thu}).
$\square$

%***********************************************************************
\subsection{Baum-Connes morphism and assembly map}
%***********************************************************************
In this section we review main aspects of Baum-Connes theory connecting
analytical and topological K-theory of discrete groups. For an excellent
introductory material we recommend  the survey of Higson (\cite{Hig}).

\medskip\noindent
\underline{Baum-Connes morphism for discrete groups.}
Let $G$ be a countable group. The analytical part
of the Baum-Connes morphism involves algebraic $K$-groups 
of {\it reduced $C^*$-algebra $C^*_{red}G$} of group $G$.
One can think of $C^*_{red}G$ as $C^*$-algebra generated
by representation of $G$ on the Hilbert space $l^2(G)$,
see Higson (\cite{Hig}). The $K$-theory of $C^*$-algebras
is reviewed in the next section. The geometrical-topological 
part of the Baum-Connes morphism involves the $K$-homology
of a {\it classifying space $BG$} for group $G$. In particular,
when $G$ has no torsion elements (which is always true in
the case of Kleinian groups), then $BG=K(G,1)$
the Eilenberg-Mac Lane space of group $G$. The {\it
Baum-Connes morphism} is a mapping
%**************************************************************
\begin{equation}\label{eq24*}
\mu_i^G: K^G_i(BG)\longrightarrow K_i(C^*_{red}G), \qquad i=0,1.
\end{equation}
%******************************************************************
(Roughly speaking, injectivity of $\mu_i^G$ allows to obtain topological
invariants of $BG$ by  looking at $C^*$-algebras associated to $G$. 
Similarly, the surjectivity of $\mu_i^G$ gives an opportunity to
classify certain $C^*$-algebras using topological data carried by the
space $BG$.)

\medskip\noindent
\underline{Assembly (index) map.}
The pointwise correspondence $\mu^G_i$ between elements of $K$-groups
to both sides of morphism (\ref{eq24}) is called an {\it assembly
map}. There exists no universal recipe of how to construct such a
map for the concrete countable groups. However, for many topologically
significant groups (e.g. amenable, hyperbolic and fundamental groups
of Haken 3-manifolds) it was done explicitly. We conclude this section by the following lemma, which 
illustrates how theory of geodesic laminations and Kleinian groups 
helps to construct the assembly map
for the fundamental group of mapping tori with a pseudo-Anosov 
monodromy.  
%*************************************************************
\begin{lem}\label{lm7}
Let $\pi_1M$ be fundamental group of the surface bundle with pseudo-Anosov
monodromy $\varphi$, such that its representation as Kleinian group $G$
is given by a $\varphi$-invariant minimal geodesic lamination $\Lambda_G=BG$. 
Then there exists an injective  (pre-) assembly map 
%**************************************************************
\begin{equation}\label{eq24}
\mu^G: \Lambda_G\longrightarrow C^*_{SAF}G,
\end{equation}
%******************************************************************
where $C^*_{SAF}$ is the category of $AF$ $C^*$-algebras of stationary
type (see Section 4.3).
\end{lem}
%****************************************************************
{\it Proof.} The statement follows our results in (\cite{Nik}). The proof
is based on combinatorics of $AF$ $C^*$-algebras (Bratteli diagrams)
and  consists in the identification of infinite paths of Bratteli
diagrams with leaves of minimal lamination $\Lambda_G$. The main
technical idea comes from symbolic dynamics and uses the
Koebe-Morse coding of geodesic lines. 
$\square$

%***********************************************************************
\subsection{K-theory of AF $C^*$-algebras}
%***********************************************************************
This section is reserved for the basic facts of $K$-theory of the 
$AF$ $C^*$-algebras. Most references  can be found in the monograph
of Rordam, Lautstsen (\cite{}). Stationary $K_0$-groups and their
classification in terms of the subshifts of finite type are discussed
by Effros in (\cite{E}). Hyperbolic representation and parametrization
of dimension groups by continued fractions can be found in (\cite{Nik}).

%****************************************************************
\medskip\noindent
\underline{Dimension groups.}
%*****************************************************************
Let $A$ be a unital $C^*$-algebra and $V(A)$ 
 a matrix $C^*$-algebra with entries in $A$.
Projections $p,q\in V(A)$ are equivalent if there exists a partial
isometry $u$ such that $p=u^*u$ and $q=uu^*$. The equivalence
class of projection $p$ is denoted by $[p]$.
Equivalence classes of orthogonal projections can be made to
a semigroup by putting $[p]+[q]=[p+q]$. The Grothendieck
completion of this semigroup to an abelian group is called
a {\it $K_0$-group of algebra $A$}. 
Functor $A\to K_0(A)$ maps a category of unital
$C^*$-algebras into the category of abelian groups so that
the positive elements $A^+\subset A$ correspond to a ``positive
cone'' $K_0^+\subset K_0(A)$ and the unit element $1\in A$
corresponds to an ``order unit'' $[1]\in K_0(A)$.
The ordered abelian group $(K_0,K_0^+,[1])$ with an order
unit is called a {\it dimension (Elliott) group} of $C^*$-algebra $A$.

%*********************************************************************
\medskip\noindent
\underline{Representation of dimension groups by geodesic lines.}
%**********************************************************************
For $k=1,2,\dots $, let us consider dimension groups of rank $2k$, i.e. dimension groups
based on the lattice ${\Bbb Z}^{2k}$.  Such dimension groups can be 
represented (and classified) by geodesic lines on the Riemann surface
of genus $k$. Let us show first that any simple (with no self-intersection)
non-periodic geodesic gives rise to a dimension group of rank $2k$. Fix
a Riemann surface $S$ of genus $k$ together with a point $p\in S$ and
simple geodesic $\gamma$ through $p$. Consider the set
%*************************************************************
\begin{equation}
Sp ~(\gamma)=\{\gamma_0,\gamma_1,\gamma_2,\dots\},
\end{equation}
%************************************************************   
of periodic geodesics $\gamma_i$ based in $p$, which monotonically 
approximate $\gamma$ in terms of ``length'' and ``direction''. 
The set $Sp ~(\gamma)$  is known as {\it spectrum} of $\gamma$
and is defined uniquely upon $\gamma$. Let $H_1(S;{\Bbb Z})={\Bbb Z}^{2k}$
be the integral homology group of surface $S$. Since each $\gamma_i$
is a 1-cycle, there is an injective map $f:Sp~(\gamma)\to 
H_1(S;{\Bbb Z})$, which relates every closed geodesic its
homology class. Note that $f(\gamma_i)=p_i\in {\Bbb Z}^{2k}$
is ``prime'' in the sense that it is not an integer multiple of some
other point of lattice ${\Bbb Z}^{2k}$. Denote by $Sp_f(\gamma)$ 
the image of $Sp~(\gamma)$ under mapping $f$. Finally, let 
$SL(2k,{\Bbb Z})$ be the group of $2k\times 2k$ integral matrices
of determinant $1$ and $SL(2k,{\Bbb Z}^+)$ its semigroup
consisting of matrices with strictly positive entries. 
It is not hard to show, that in appropriate basis in 
$H_1(S;{\Bbb Z})$ the following is true: (i) the coordinates
of vectors $p_i$ are non-negative; (ii) there exists a matrix
$A_i\in SL(2k,{\Bbb Z}^+)$ such that $p_i=A_i(p_{i-1})$ for
any pair of vectors  $p_{i-1},p_i$ in $Sp_f(\gamma)$. The dimension group
$({\Bbb Z}^{2k},({\Bbb Z}^{2k})^+,[1])$ defined as inductive limit
of simplicial dimension groups:
%*********************************************************************
\begin{equation}\label{d.g.}
\Z^{2k}\buildrel\rm A_1\over\longrightarrow \Z^{2k}
\buildrel\rm A_2\over\longrightarrow
\Z^{2k}\buildrel\rm A_3\over\longrightarrow \dots,
\end{equation}
%*******************************************************************
is called {\it  associated} to geodesic $\gamma$.  In fact, every
dimension group of rank $2k$ can be obtained in such a  way (we omit
the proof of this fact here).

%*******************************************************************************
\medskip\noindent
\underline{Parametrization of dimension groups by continued fractions.}
%*******************************************************************************
It is known that dimension groups of rank $2$ are parametrized by the infinite
continued fraction converging to a ``slope''
of non-closed geodesic on 2-torus, which characterizes such groups, see e.g. 
(\cite{E}). Similar description is available for the dimension groups of rank
$k\ge 2$. Let $(S,\gamma)$ be the pair consisting of the Riemann surface $S$
and non-periodic geodesic $\gamma$ which represents dimension group
$({\Bbb Z}^{2k},({\Bbb Z}^{2k})^+,[1])$ as described in previous paragraph.  
Intuitively, every homotopy class of non-periodic simple lines has
a ``slope'' on $S$ which is given by a real irrational number $\theta$.
(Existence of such a number was suggested by Andr\'e Weil in 1933.)
This number was defined in (\cite{Nik}). We assume $G=\pi_1(S)$ is the principal 
congruence subgroup $\Gamma(n)$ of group $PSL_2({\Bbb Z})$.

Recall that the geodesic lines on surface $S={\Bbb H}/G$ fall into 
(i) quasi-ergodic, (ii) periodic and (iii) quasi-periodic classes  
(Artin, Hedlund, Hopf). Component (i) consists of geodesics (with self-intersections) 
which approximate any geodesic on $S$ issued from any point in any direction.
Such geodesics are ``typical'' in the sense of measure at the boundary $\partial {\Bbb H}$.
To the contrary, geodesics of class (ii) and (iii) aren't quasi-ergodic
and can be thought of as having specific ``direction'' on $S$. Non-quasi-ergodic
geodesics have measure zero and cardinality of continuum at $\partial {\Bbb H}$. 
Class (ii) is countable and class (iii) is an uncountable set of simple (without
self-intersections) geodesics and they are a substitute of {\it global} straight 
lines on $S$. Class (iii) will be basic for our definition of ``slope'' since it 
consists of simple self-approximating geodesics with definite direction on $S$.

Let $\Delta$ be the fundamental domain of $G$. Then every geodesic $\gamma$
of class (iii) intersects $\Delta$ infinitely often, and let us enumerate
the corresponding arcs on intersection by $\{[\gamma_0],[\gamma_1],\dots\}$
in the order they cut $\Delta$. Since each $[\gamma_i]$ is the arc
of the geodesic half-circle on $\Bbb H$, they can be moved to each other
by transformations from $G$. Put then $g_i\in G$ be such that $[\gamma_i]=
g_i([\gamma_{i-1}])$. It can be shown that $g_i\in \Gamma^+(n)$, where
$\Gamma^+(n)$ is the principal congruence semi-group consisting of 
$2\times 2$ matrices with positive integer entries.  
Finally, suppose that $g_i=\left(\small\matrix{a_0^{(i)} & 1\cr 1 & 0}\right)\dots
\left(\small\matrix{a_n^{(i)} & 1\cr 1 & 0}\right)$ is the Minkowski
decomposition of matrix $g_i\in G$. Then regular continued fraction

\bigskip\noindent
%***************************************************************************
\begin{equation}\label{e12*}
\theta=a_0^{(0)}+{1\over\displaystyle a_1^{(0)}+
{\strut 1\over\displaystyle a_2^{(0)}\displaystyle +\dots}},
\end{equation}
%*****************************************************************

\medskip\noindent
is called {\it associated} to the dimension group defined by pair $(G,\gamma)$.
By analogy with the case of noncommutative torus, real number $\theta$
is called a {\it rotation number}. Two dimension groups are order-isomorphic 
if and only if their continued
fractions coincide, except possibly in the finite number of terms.
In other words, the corresponding rotation numbers are modular equivalent:
$\theta'={a\theta+b\over c\theta +d},
\qquad a,b,c,d\in\Z,\qquad ad-bc=\pm 1$.

%******************************************************************
\medskip\noindent
\underline{Stationary dimension groups and their classification.}
%*******************************************************************
Dimension group $({\Bbb Z}^{2k},$\linebreak $({\Bbb Z}^{2k})^+,[1])$ is
called {\it stationary} if in formula (\ref{d.g.})  $A_i=Const=A$.
If continued fraction (\ref{e12*}) is periodic, then the corresponding
dimension group is stationary. Indeed, in this case there exists
a transformation $g\in G$ such that $g_i\dots g_{i+k}=~Const=g$,
where $k$ is the minimal period of continued fraction. It isn't
hard to see that $g$ is a hyperbolic transformation with two
fixed points at the boundary $\partial {\Bbb H}$. The geodesic
$\gamma$ through these points is invariant under the action of
$g$ and coincides with the geodesic corresponding to the stationary
dimension group generated by the automorphism $A$. Conversely,
each stationary dimension group gives rise to a geodesic, whose
slope (rotation number) unfolds into a periodic continued fraction.

Since periodic continued fractions converge to quadratic irrationalities,
one can classify stationary dimension groups using the arithmetic
of the real quadratic number fields. The corresponding classification
scheme was suggested in Section 1.

%**************************************************************************

%**********************************************************

\begin{thebibliography}{100}



\bibitem{Bia}
L.~Bianchi, Sui gruppi di sostituzioni lineari con coefficienti
appartenenti a corpi quadratici immaginari, Math. Ann. 40 (1892),
332-412.



\bibitem{Bor}
A.~Borel, Commensurability classes and volumes of hyperbolic 3-manifolds,
Ann. Scuola Norm. Sup. Pisa 8 (1981), 1-33.


\bibitem{BMR}
B.~H.~Bowditch, C.~Maclachlan and A.~W.~Reid, 
Arithmetic hyperbolic surface bundles, Math. Ann. 302 (1995), 31-60.

\bibitem{Boy}
D.~W.~Boyd, Mahler's measure and invariants of hyperbolic manifolds. 
Number theory for the millennium, I (Urbana, IL, 2000), 127-143.

\bibitem{Bro}
J.~F.~Brock, The Weil-Petersson metric and volumes of 3-dimensional
hyperbolic convex cores, J. Amer. Math. Soc. (to appear)


\bibitem{E}
E.~G.~Effros, Dimensions and $C^*$-Algebras, Conf. Board of the Math.
Sciences No.46, AMS (1981).

\bibitem{Han}
D.~Handelman, Positive matrices and dimension groups affiliated
to $C^*$-algebras and topological Markov chains, J. Operator
Theory 6 (1981), 55-74.


\bibitem{H}
E.~Hecke, Vorlesungen \"uber die Theorie der Algebraischen Zahlen,
Chelsea, N.Y., 1948.

\bibitem{Hel}
H.~Helling, Bestimmung der Kommensurabilit\"atklasse der Hilbertschen
Modulgruppe, Math. Zeitschr 92 (1966), 269-280.


\bibitem{HPS}
R.~H.~Herman, I.~F.~Putnam and C.~F.~Skau, Ordered Bratteli diagrams,
dimension groups and topological dynamics, Internat. J. Math. 3 (1992), 
827-864.



\bibitem{Hig}
N.~Higson, The Baum-Connes conjecture, Doc. Math. J. DMV,
Extra volume ICM 1998, II, 637-646. 

\bibitem{Hum}
G.~Humbert, Sur la mesure des classes d'Hermite de discriminante
donn\'e dans un corps quadratique imaginaire, et sur certaines
volumes non euclidiens, C.~R.~Acad. Sci. Paris 169 (1919), 448-454.

\bibitem{Kas}
R.~M.~Kashaev, A link invariant from quantum dilogarithm. Modern Phys. 
Lett. A 10 (1995), no. 19, 1409-1418. 


\bibitem{McM}
C.~T.~McMullen, Renormalization and 3-Manifolds which Fiber over the
Circle, Annals of Mathematical Studies 142, Princeton, N.~J., 1996.


\bibitem{NeZ}
W.~D.~Neumann and D.~Zagier, Volumes of hyperbolic three-manifolds,
Topology 24 (1985), 307-332.


\bibitem{Nik}
I.~Nikolaev, Geodesic laminations and noncommutative geometry, math.GT/0209155,
26p.



\bibitem{RLL}
M.~R\o rdam, F.~Larsen and N.~Laustsen, An introduction to $K$-theory 
for $C^*$-algebras. London Mathematical Society Student Texts, 49. 
Cambridge University Press, Cambridge, 2000.


\bibitem{Thu}
W.~P.~Thurston, Three dimensional manifolds, Kleinian groups
and hyperbolic geometry, Bull. Amer. Math. Soc. 6 (1982), 357-381.

\bibitem{T}
W.~P.~Thurston, The Geometry and Topology of Three-Manifolds,
MSRI 1997, electronic edition of 1980 Princeton Univ. notes,
available at {\sf http://www.msri.org/gt3m/}


\bibitem{Wee}
J.~Weeks, SNAPPEA, available {\sf www northnet.org/weeks/index/SnapPea.html}

\end{thebibliography}
\end{document}